# Miniaturization-Oriented Design of Spline-Parameterized UWB Antenna for In-Door Positioning Applications

Adrian Bekasiewicz[1][0000-0003-0244-541X], Tom Dhaene[2][0000-0003-2899-4636], Ivo Couckuyt[2][0000-0002-9524-4205], and Jacek Litka[1][0000-0002-3538-4118]

[1] Faculty of Electronics, Telecommunications and Informatics, Gdansk University of Technology, Narutowicza 11/12, 80-233 Gdansk, Poland
[2] Department of Information Technology (INTEC), IDLab, Ghent University-imec, iGent, Technologiepark-Zwijnaarde 126, 9052 Ghent, Belgium
`bekasiewicz@ru.is`

**Abstract.** Design of ultra-wideband antennas for in-door localization applications is a challenging task that involves development of geometry that ensures appropriate balance between the size and performance. In this work, a topologically-flexible monopole has been generated using a stratified framework which embeds a gradient-based trust-region (TR) optimization algorithm in a meta-loop that gradually increases structure dimensionality. The optimization has been performed using a composite objective function that maintains acceptable size/performance trade-off. The final design features a reflection below –10 dB within the UWB spectrum and a small footprint of only 182 mm$^2$. The considered method has been benchmarked against a standard TR-based routine executed directly on a multi-dimensional representation of the antenna model.

## 1    Introduction

Ultra-wideband (UWB) technology is a promising solution for the development of modern real-time localization services dedicated to in-door environments. Its advantage—compared to other radio-frequency methods—includes pulse operation over a broad frequency range which makes it less susceptible to interferences [1], [2]. Quality-of-service offered by positioning systems is a subject to availability of antennas which comply with the regulations that mandate the access to wireless medium while ensuring high performance [3]. Apart from the electrical- and field-related requirements, applicability of UWB radiators in mobile terminals is also affected by their footprint [4]. Therefore, accounting for the size/performance trade-off is important when design of radiators for in-door positioning is considered [1], [4].

Conventional approaches to design of compact antennas involve experience-driven determination of the topology followed by its trial-and-error-based modifications in hope of achieving improved performance and/or area reduction [5], [6]. This inherently cognitive process is affected by the engineering bias as designers, understandably, lean towards geometries and modifications they are familiar with [4], [7]. Although the structure development should be followed by rigorous optimization-based tuning, this step is often neglected in favor of parametric studies (aided by visual-inspection of performance changes) oriented towards achieving the satisfactory performance [8]. The main reasoning behind this laborious and prone-to-failure procedure is that accurate evaluation of complex antenna performance can only be performed using computationally expensive electromagnetic (EM) simulations. At the same time, a large number of EM simulations required to converge challenges the applicability of conventional algorithms. The undesirable consequence of the outlined scheme is that the antenna design is often governed by a mix of past experiences and limited number of observations which, at best, lead to acceptable solutions rather than the ones that offers the best balance between the multiple requirements [5], [7], [8].

The effects of engineering-bias on the design process can be mitigated by shifting the design paradigm from the fixed-topology to flexible models where the final geometry is determined by numerical methods. This can be achieved using so-called dummy EM models which represent the antenna as a set of interconnected points, or in the form of a binary matrix that govern the structure shape [9], [10], [11]. The bottleneck of generic models is that they require overwhelmingly large number of dimensions to support diverse geometries and hence the cost of their EM-driven optimization is numerically prohibitive. The challenges pertinent to unacceptable optimization cost can be mitigated using trust-region (TR) methods. TR algorithms exploit the promising regions of the search space based on evaluations of data-efficient models that approximate EM simulations [12], [13]. Despite being numerically-efficient, local nature of TR-based optimization might limit its applicability for solving complex and multi-modal problems such as automatic tuning of topology oriented towards ensuring a balance between the size and performance.

In this work, a framework for stratified design of topologically-flexible antenna has been considered. The method embeds the TR-based optimization into a meta-loop that enables a gradual increase of the structure dimensionality so as to first identify and then exploit the promising region of the search space. The design process is governed by a composite objective function that balances the trade-off between the antenna size and its performance. The method has been used for the development of a spline-parameterized UWB monopole dedicated to in-door positioning applications. The optimized structure features a reflection at most –10 dB within 3.1 GHz to 10.6 GHz bandwidth and a footprint of only 182 mm$^2$. The method has been benchmarked against a standard TR-based optimization.

## 2      Spline-Parameterized Monopole

Optimization of the radiator in a dimensionally-flexible setup is subject to availability of a suitable EM simulation model. Here, the UWB antenna of Fig. 1 has been considered [14]. The structure is designed on a substrate with permittivity/thickness of 3.38/0.813 mm. It features a driven element in the form of a radiator fed through a microstrip line and a ground plane with an L-shaped extension. The antenna EM sim-ulation model is implemented in CST Microwave Studio and evaluated using its time-domain solver [15]. The structure is discretized using 400,000 hexahedral mesh cells (on average), whereas its typical simulation time amounts to 160s.

In order to support dimensionality adjustment, the radiator and ground plane are represented in the form of spline-based curves defined using variable number of points (so-called knots) [11]. The vector of design parameters is $x = [x_c\ Y \cdot x_g\ S \cdot x_r]^T$, where $x_c = [X\ l_f\ l_1\ l_{2r}\ w_1\ o_r]^T$ represents topology-specific variables; $x_g = [x_{g.1}\ \ldots\ x_{g.l}\ \ldots\ x_{g.L}]^T$ and $x_r = [x_{r.1}\ \ldots\ x_{r.l}\ \ldots\ x_{r.L}]^T$ are the coordinates of the knots ($l = 1, \ldots, L$) spanned equidistantly between the antenna edges (for the ground plane) and along the $2\pi$ azimuth angles defined in the cylindrical coordinates (for the radiator; cf. Fig. 1). Maintaining equidistant distribution of points prevents self-intersections of the generated splines which is crucial to ensure feasible topologies. The parameters $l_2 = (X - w_1)l_{2r}$, $Y = l_1 + w_1$, $l_{fr} = \min(X, Y - l_f)/2$, $S = \min(X - o_r, Y - l_f)/2$, and $o = 0.5X + o_r$ are relative, whereas $w_f = 1.8$ is set to ensure 50 Ohm input impedance. Note that $l_{2r}$ is dimensionless and all other variables are in mm. The antenna model is considered feasible within the following lower and upper bounds $l_b = [6\ 4\ 10\ 0.05\ 0.5\ -1\ 0.2\mathbf{1}^L\ 0.1\mathbf{1}^L]^T$ and $u_b = [30\ 15\ 30\ 1\ 2.5\ 1\ 0.8\mathbf{1}^L\ \mathbf{1}^L]^T$, where $\mathbf{1}^L$ is the $L$-dimensional vector of ones. The vector $x$ is represented using $2L + 6$ dimensions. For each EM simulation, the antenna is generated dynamically using the specified number of points [15].

## 3  Design Methodology

### 3.1  Problem Formulation

Let $R(x) = R(x, f)$ be the response of the antenna obtained over a frequency sweep $f$ for the vector of input parameters $x$. The optimization task is given as:

$$x^* = \mathrm{argmin}\left(U\left(R(x)\right)\right) \qquad (1)$$

where $x^*$ is the optimal design to be found and $U(x) = U(R_f(x))$ is a scalar objective function. Direct solving of (1) is numerically impractical when multi-parameter structures are considered. Instead, the problem can be solved using a gradient-based algorithm embedded within the TR framework [12].

### 3.2  TR-Based Optimization Algorithm

The considered optimization engine generates a series of approximations, $i = 0, 1, 2, \ldots$, to the original problem by solving:

$$x^{(i+1)} = \arg\min_{\|x - x^{(i)}\| \leq \delta^{(i)}} \left(U\left(G^{(i)}(x)\right)\right) \qquad (2)$$

where $G^{(i)} = R(x^{(i)}) + J(x^{(i)})(x - x^{(i)})$ is a first-order Taylor expansion model, whereas $J$ is a Jacobian generated around the $x^{(i)}$ using a large-step finite differences [2], [12]. The TR radius, that determines the region of $G^{(i)}$ model validity (around $x^{(i)}$), is ad-justed based on the gain coefficient $\rho = [U(\mathbf{R}(\mathbf{x}^{(i+1)})) - U(\mathbf{R}(\mathbf{x}^{(i)}))]/[U(\mathbf{G}^{(i)}(\mathbf{x}^{(i+1)})) - U(\mathbf{G}^{(i)}(\mathbf{x}^{(i)}))]$ which expresses the expected versus obtained change of the objective function. The radius is initialized as $\delta^{(0)} = 1$ and then adjusted according to gain factor as $\delta^{(i+1)} = 2\delta^{(i)}$ when $\rho > 0.75$, and $\delta^{(i+1)} = \delta^{(i)}/3$ when $\rho < 0.75$. The algorithm is terminated when $\|\mathbf{x}^{(i+1)} - \mathbf{x}^{(i)}\| \leq \varepsilon$, or $\delta^{(i+1)} \leq \varepsilon$, where $\varepsilon = 10^{-2}$ ($\|\cdot\|$ represents the Euclidean norm). Note that the method is data-efficient as it requires only $D + 1$ EM simulations ($D$ is the problem dimensionality) for construction of $\mathbf{G}^{(i)}$ per successful iteration, i.e., when $\rho > 0$. Additional EM model evaluation is required for each unsuccessful step. For more comprehensive discussion on the TR-based optimization, see [2], [12], [14].

### 3.3    Stratified Design Framework

The TR algorithm is embedded within a stratified design framework which gradually increases the number of antenna parameters [14]. The goal of the process is to first identify (using low number of parameters) and then exploit (by increasing problem dimensionality) the promising regions of the search space. The method involves a series of TR-based optimizations of the antenna of Section 2 being represented using a specified number of variables. Let $\mathbf{L} = [L^{(0)} \ldots L^{(j)} \ldots L^{(J)}]^T$ ($j = 0, 1, \ldots, J$) be a vector that defines the number of spline parameters used to construct antenna ground plane and radiator in consecutive meta-iterations. The algorithm is initialized using $L^{(0)}$, whereas at the beginning of each consecutive iteration ($j > 0$) the vectors $\mathbf{x}_g^{(j+1)}$ and $\mathbf{x}_r^{(j+1)}$ required for construction of $\mathbf{x}_0^{(j+1)}$ are obtained though interpolation of their counterparts extracted from $\mathbf{x}_{\text{opt}}^{(j)}$ found as a result of previous meta-iteration. The considered stratified design framework can be summarized as follows:

1. Specify $L$, set $j = 0$ and define $\mathbf{x}_0^{(j)}$;
2. Optimize $\mathbf{x}_0^{(j)}$ using algorithm of Section 3.2 to obtain $\mathbf{x}_{\text{opt}}^{(j)} = \mathbf{x}^*$;
3. If $j = J$ or $U(\mathbf{x}_{\text{opt}}^{(j+1)}) \leq U(\mathbf{x}_{\text{opt}}^{(j)})$, set $\mathbf{x}^* = \mathbf{x}_{\text{opt}}^{(j)}$ and END; otherwise obtain $\mathbf{x}_g^{(j+1)}$, $\mathbf{x}_r^{(j+1)}$ by $L^{(j+1)}$-point interpolation of $\mathbf{x}_g^{(j)}$, $\mathbf{x}_r^{(j)}$ extracted from $\mathbf{x}_{\text{opt}}^{(j)}$ and go to Step 4;
4. Extract $\mathbf{x}_c^{(j)}$ from $\mathbf{x}_{\text{opt}}^{(j)}$ and define $\mathbf{x}_0^{(j+1)} = [\mathbf{x}_c^{(j)} \mathbf{x}_g^{(j+1)} \mathbf{x}_r^{(j+1)}]^T$, set $j = j + 1$ and go to Step 2.

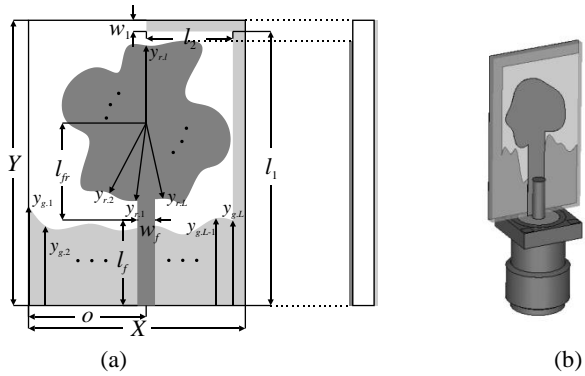

**Fig. 1.** A spline-parameterized monopole: (a) geometry with highlight on the design parameters and (b) visualization of the structure. Note that $y_{g.l} = Y \cdot x_{g.l}$ and $y_{r.l} = S \cdot x_{r.l}$ ($l = 1, \ldots, L$).

It should be noted that the number of $L^{(j)}$ points affects the optimization cost as the number of EM simulations required for construction of the linear model is proportional to problem dimensionality (cf. Section 2.3). On the other hand, embedding the design task in stratified framework mitigates the risk of getting stuck in a poor local optimum.

### 3.4   Objective Function

The objective function considered for optimization is given as (see Fig. 2 for conceptual illustration) [14]:

$$U(x, \alpha) = \begin{cases} U_1(x), \text{ when } \alpha = 0 \\ A(x) + \gamma_1 \max(U_1(x)/S_1, 0), \text{ when } \alpha = 1 \\ \max(S(x)) + \gamma_2 \max((A(x) - A_1^{(i)})/A_1^{(i)}, 0), \text{ when } \alpha = 2 \end{cases} \quad (3)$$

where $S(x) = |R(x)|_{f_L \leq f \leq f_H}$ represents the structure reflection response (in dB) within $f_L$ = 3.1 GHz to $f_H$ = 10.6 GHz, threshold $S_1$ = –10 dB, $A(x) = S \cdot Y$ denotes the antenna footprint (cf. Section 2), whereas $\gamma_1$ = 1000 and $\gamma_2$ = 500 are the scaling coefficients; $U_1(x) = \max(S(x)) - S_1$. The thresholds on the size $A_1^{(i)}$ are recorded when the function selector $\alpha$ is triggered.

The reasoning behind the use of composite objective (3) is to ensure the acceptable balance between antenna size and its electrical performance. When $\alpha = 0$ the optimization is oriented only towards minimization of the in-band reflection. Once attained (i.e., $U(x, \alpha = 0) < 0$), the algorithm records $A_1^{(i)}$ sets $\alpha = 1$ and performs explicit miniaturization of $A$ with penalty on $S_{11}$. When $\max(S(x)) > S_2$ (here, $S_2 = -9.5$ dB; cf. Fig. 2), the mode is set to $\alpha = 2$, $A_1^{(i)}$ is stored, and the explicit minimization of reflection with penalty on size is performed until $\max(S(x)) < S_1$, when $\alpha = 1$ is selected again.

## 4   Numerical Results

The antenna is optimized using the framework of Section 3. The goal of the design process is to miniaturize the structure while maintaining its in-band reflection at the level of around –10 dB (controlled by $S_1$, $S_2$ thresholds). The vector of spline-knots for consecutive meta-iterations and the initial design are $L$ = [1 8 16 24 32]$^T$, and $x_0^{(0)}$ = [10 6 16 0.8 1 0 0.35 0.6]$^T$, respectively. The final design $x^* = x_{opt}^{(2)}$ = [10.99 4.87 15.37 1 1.53 –0.31 0.45 0.49 0.24 0.2 0.47 0.2 0.21 0.35 0.52 0.49 0.78 0.39 0.89 0.68 0.59 0.6]$^T$ is obtained after three meta-iterations. The algorithm termination is triggered at $j = 2$, due to lack of the objective function improvement.

The optimized design is characterized by reflection below –10 dB within the entire bandwidth and the dimensions of 10.99 mm × 16.91 mm. Although the final solution is larger by 18 mm$^2$ compared to the starting point (the footprint of 182.4 mm$^2$ at $x^*$ vs. 170 mm$^2$ at $x_0^{(0)}$), its in-band reflection is –5.7 dB lower—w.r.t. to maximum in-band reflection of –4.2 dB at $x_0^{(0)}$—which demonstrates flexibility of the discussed framework and composite objective function in terms of balancing the electrical performance and physical dimensions of the radiator. The convergence curves for the TR algorithm at the first two meta-steps, as well as the reflection responses at the selected designs are shown in Fig. 3. The geometry of the optimized structure and its radiation-pattern characteristics obtained at 4 GHz and 6.5 GHz (in xy-plane) are given in Fig. 4. Note that the considered test frequencies are centered w.r.t. channels 2 and 5 of the UWB band [3]. The radiation patterns are fairly omnidirectional. Their distortion around the –90° direction is due to implementation of the L-shaped ground plane extension which supports miniaturization, but also acts as a signal reflector (cf. Fig. 1). Nonetheless, owing to small size the radiator is suitable for wearable applications.

The considered framework has been compared against a conventional TR-based algorithm in terms of cost and performance. The benchmark method has been used for direct optimization of the antenna represented by a 38-dimensional vector of parameters. The process was controlled only by the objective function (3) with $\alpha = 1$. For fair comparison, the same initial design $x_0^{(0)}$ has been set for both algorithms. The results gathered in Table 1 indicate that, although the benchmark algorithm obtained smaller design at a lower cost, it failed at fulfilling the performance-related specifications. From this perspective, the considered stratified design approach offers improved balance between the antenna size and its electrical properties.

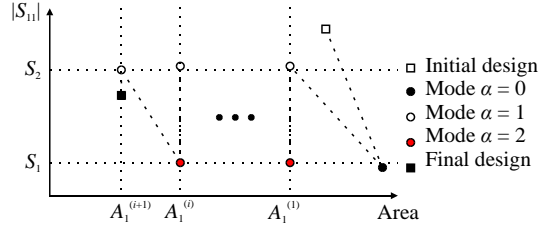

**Fig. 2.** TR-based optimization of the antenna using a composite objective function that balances the electrical- and size-related specifications.

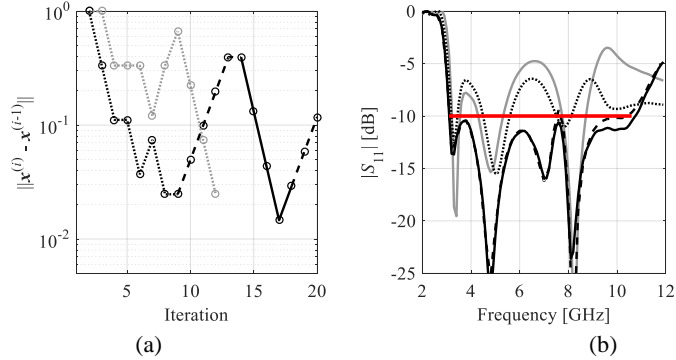

**Fig. 3.** A spline-based monopole: (a) convergence plots for the first 20 TR-based iterations at the meta-steps $j = 0$ (gray) and $j = 1$ with $\alpha = 0$ (···), $\alpha = 1$ (– –) and $\alpha = 2$ (–), as well as (b) responses at $x_0^{(0)}$ (gray) and $x_{opt}^{(0)}$ (···), $x_{opt}^{(1)}$ (—), and $x_{opt}^{(2)}$ (– –) designs. Note that the performance $x_{opt}^{(2)}$ is slightly worse w.r.t. $x_{opt}^{(1)}$ due to interpolation between meta-iterations.

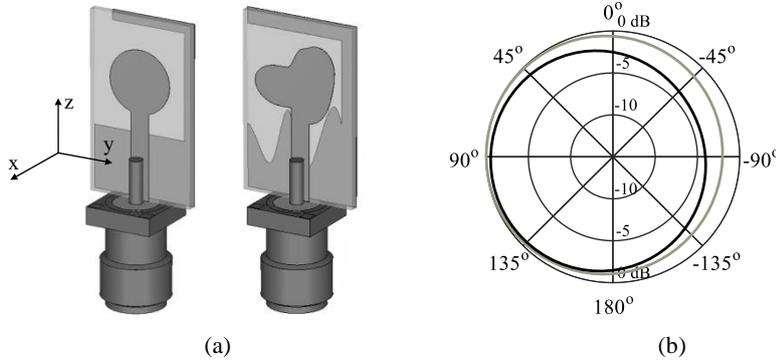

**Fig. 3.** A spline-based antenna: (a) in-scale comparison of geometries at $x_0^{(0)}$ (left) and $x^*$ (right) designs, as well as (b) xy-plane radiation patterns at 4 GHz (gray) and 6.5 GHz (black).

**Table 1.** Benchmark of the stratified design optimization framework.

| Method | Meta-step – cost [$R$] | | | Total cost [$R$] | Size [mm$^2$] | max($S(x)$) [dB] |
|---|---|---|---|---|---|---|
| | 1 | 2 | 3 | | | |
| TR with $\alpha = 1$ | 241 | – | – | 241 | 171.6 | –8.5* |
| This work | 76 | 491 | 39 | 606 | 182.4 | –10 |

\* The optimized design violates the performance specifications resulting from $S_1$

## 5 Conclusion

In this work, a framework dedicated to design of topologically flexible antenna structures has been discussed. The method involves a series of TR-based optimizations embedded within a meta-loop that gradually increases problem dimensionality. The TR-based design is governed by a composite objective function that balances the size/performance trade-off. The final design is characterized by a reflection of –10 dB within the UWB bandwidth (i.e., from 3.1 GHz to 10.6 GHz), as well as a footprint of only 182 mm$^2$. The far-field responses evaluated at 4 GHz and 6.5 GHz frequencies (that represent channel 2 and 5 of the UWB spectrum) indicate that the antenna features a fairly omnidirectional radiation patterns. Small dimensions and high performance make the structure of potential use for wearable in-door positioning devices. The benchmark of the considered framework against a standard TR-based optimization demonstrates that, despite higher computational cost, the method offers improved balance between the performance figures and dimensions. Future work will focus on application of the antenna prototypes as components of a real-world in-door positioning system.

**Acknowledgments.** This work was supported in part by the National Science Centre of Poland Grants 2020/37/B/ST7/01448 and 2021/43/B/ST7/01856, National Centre for Research and Development Grant NOR/POLNOR/HAPADS/0049/2019-00, and Gdansk University of Technology (Excellence Initiative - Research University) Grant 16/2023/IDUB/IV.2/EUROPIUM.

**Disclosure of Interests.** The authors declare no conflicts of interest.